\documentclass[12pt]{article}

\makeatletter
\newfont{\footsc}{cmcsc10 at 8truept}
\newfont{\footbf}{cmbx10 at 8truept}
\newfont{\footrm}{cmr10 at 10truept}
\renewcommand{\ps@plain}{%
\renewcommand{\@oddfoot}{\footsc the electronic journal of combinatorics
  {\footbf 8} (2001), \#N9\hfil\footrm\thepage}}
\makeatother 

\title{A LOWER BOUND FOR THE NUMBER OF EDGES IN A GRAPH CONTAINING NO TWO CYCLES
OF THE SAME LENGTH }

\author{ Chunhui Lai \thanks{ Project Supported by NSF of Fujian(A96026),
Science and Technology Project of Fujian(K20105) and Fujian
Provincial Training Foundation
for "Bai-Quan-Wan Talents Engineering".}\\
\small Dept. of Math., Zhangzhou Teachers College, \\[-0.8ex]
\small Zhangzhou, Fujian 363000, P. R. of CHINA. \\[-0.8ex]
\small \texttt{ zjlaichu@public.zzptt.fj.cn }}

\date{\small Submitted: November 3, 2000;  Accepted: October 20, 2001.\\
\small MR Subject Classifications: 05C38, 05C35}

\begin{document}
\maketitle

\begin{abstract}
In 1975, P. Erd\"{o}s proposed the problem of determining the
maximum number $f(n)$ of edges in a graph of $n$ vertices in which
any two cycles are of different
 lengths. In this paper, it is proved that $$f(n)\geq n+32t-1$$ for
$t=27720r+169 \,\ (r\geq 1)$
 and $n\geq\frac{6911}{16}t^{2}+\frac{514441}{8}t-\frac{3309665}{16}$.
Consequently, $\liminf\sb {n \to \infty} {f(n)-n \over \sqrt n}
\geq \sqrt {2 + {2562 \over 6911}}.$
 \end{abstract}

\section{Introduction}
Let $f(n)$ be the maximum number of edges in a graph on $n$
vertices in which no two cycles have the same length. In 1975,
Erd\"{o}s raised the problem of determining $f(n)$ (see [1],
p.247, Problem 11). Shi[2] proved that
$$f(n)\geq n+
 [(\sqrt {8n-23} +1)/2]$$ for $n\geq 3$. Lai[3,4,5,6] proved that for $n\geq
(1381/9)t^2 + (26/45)t + 98/45, t= 360q
 + 7,$  $$f(n)\geq n+19t-1,$$ and for $n\geq e^{2m}(2m+3)/4$, $$f(n)<
 n-2+\sqrt {n ln (4n/(2m+3)) +2n} +log_2 (n+6).$$ Boros, Caro, F\"{u}redi and
Yuster[7]
 proved that $$f(n)\leq n+1.98\sqrt{n}(1+o(1)).$$ Let $v(G)$ denote  the number
of
 vertices, and $\epsilon(G)$ denote  the number of edges. In this paper, we
construct a
 graph $G$ having no two cycles with the same length which leads to the
following result.
  \par
\bigskip
  \noindent{\bf Theorem.} Let $ t=27720r+169 \,\ (r\geq 1)$, then $$f(n)\geq
n+32t-1$$ for
  $n\geq
  \frac{6911}{16}t^2+\frac{514441}{8}t-\frac{3309665}{16}$.

 \section{Proof of Theorem}
  {\bf Proof.} Let $t=27720r+169,r\geq 1,$
  $ n_{t}=\frac{6911}{16}t^2+\frac{514441}{8}t-\frac{3309665}{16},$ $n\geq
n_{t}.$
  We shall show that there exists a graph $G$ on $n$ vertices with $ n+32t-1$
edges such
  that all cycles in $G$ have distinct lengths.

  Now we construct the graph $G$ which consists of a number of subgraphs: $B_i$,
  ($0\leq i\leq 21t+\frac{7t+1}{8}-58, 22t-798\leq i\leq 22t+64,
  23t-734\leq i\leq 23t+267, 24t-531\leq i\leq 24t+57, 25t-741\leq i\leq 25t+58,
  26t-740\leq i\leq 26t+57, 27t-741\leq i\leq 27t+57, 28t-741\leq i\leq 28t+52,
  29t-746\leq i\leq 29t+60, 30t-738\leq i\leq 30t+60,$ and $31t-738\leq i\leq
31t+799$).

  Now we define these $B_i$'s. These subgraphs all have a common vertex $x$,
otherwise their vertex sets are pairwise
  disjoint.
  \par
  For $\frac{7t+1}{8}\leq i\leq t-742,$ let the subgraph  $B_{19t+2i+1}$
consist of a
  cycle $$C_{19t+2i+1}= xx_i^1x_i^2...x_i^{144t+13i+1463}x$$ and eleven paths
  sharing a common vertex $x$, the other end vertices are on the
  cycle $C_{19t+2i+1}$:
  $$xx_{i,1}^1x_{i,1}^2...x_{i,1}^{(11t-1)/2}x_i^{(31t-115)/2+i}$$
  $$xx_{i,2}^1x_{i,2}^2...x_{i,2}^{(13t-1)/2}x_i^{(51t-103)/2+2i}$$
  $$xx_{i,3}^1x_{i,3}^2...x_{i,3}^{(13t-1)/2}x_i^{(71t+315)/2+3i}$$
  $$xx_{i,4}^1x_{i,4}^2...x_{i,4}^{(15t-1)/2}x_i^{(91t+313)/2+4i}$$
  $$xx_{i,5}^1x_{i,5}^2...x_{i,5}^{(15t-1)/2}x_i^{(111t+313)/2+5i}$$
  $$xx_{i,6}^1x_{i,6}^2...x_{i,6}^{(17t-1)/2}x_i^{(131t+311)/2+6i}$$
  $$xx_{i,7}^1x_{i,7}^2...x_{i,7}^{(17t-1)/2}x_i^{(151t+309)/2+7i}$$
  $$xx_{i,8}^1x_{i,8}^2...x_{i,8}^{(19t-1)/2}x_i^{(171t+297)/2+8i}$$
  $$xx_{i,9}^1x_{i,9}^2...x_{i,9}^{(19t-1)/2}x_i^{(191t+301)/2+9i}$$
  $$xx_{i,10}^1x_{i,10}^2...x_{i,10}^{(21t-1)/2}x_i^{(211t+305)/2+10i}$$
  $$xx_{i,11}^1x_{i,11}^2...x_{i,11}^{(t-571)/2}x_i^{(251t+2357)/2+11i}.$$
  \par
 From the construction, we notice that $B_{19t+2i+1}$ contains exactly
seventy-eight cycles of lengths:
 \par
 $$\begin{array}{llll}21t+i-57,& 22t+i+7,& 23t+i+210,& 24t+i,\\
  25t+i+1,& 26t+i,& 27t+i,& 28t+i-5,\\
  29t+i+3, & 30t+i+3,& 31t+i+742,& 19t+2i+1,\\
  32t+2i-51, & 32t+2i+216,& 34t+2i+209,& 34t+2i,\\
  36t+2i,&36t+2i-1,&38t+2i-6,&38t+2i-3,\\
  40t+2i+5,&40t+2i+744,& 49t+3i+1312,&42t+3i+158,\\
  43t+3i+215,& 44t+3i+209,& 45t+3i-1,& 46t+3i-1,\\
  47t+3i-7, &48t+3i-4,& 49t+3i-1,&50t+3i+746,\\
  58t+4i+1314, &53t+4i+157,& 53t+4i+215,& 55t+4i+208,\\
  55t+4i-2,& 57t+4i-7,& 57t+4i-5,& 59t+4i-2, \\
  59t+4i+740,& 68t+5i+1316,& 63t+5i+157,& 64t+5i+214,\\
  65t+5i+207,& 66t+5i-8,&67t+5i-5,&68t+5i-3,\\
  69t+5i+739, &77t+6i+1310,& 74t+6i+156,& 74t+6i+213,\\
  76t+6i+201,& 76t+6i-6,& 78t+6i-3,& 78t+6i+738,\\
  87t+7i+1309,& 84t+7i+155, &85t+7i+207,& 86t+7i+203,\\
  87t+7i-4, &88t+7i+738,& 96t+8i+1308,& 95t+8i+149, \\
  95t+8i+209,& 97t+8i+205,& 97t+8i+737,&106t+9i+1308,\\
  105t+9i+151,&106t+9i+211,& 107t+9i+946,&115t+10i+1307,\\
  116t+10i+153,&116t+10i+952, &125t+11i+1516,& 126t+11i+894,\\
  134t+12i+1522,& 144t+13i+1464.&&\end{array}$$
  \par

  Similarly, for $58\leq i\leq \frac{7t-7}{8},$ let the subgraph $B_{21t+i-57}$
consist of a
  cycle $$xy_i^1y_i^2...y_i^{126t+11i+893}x$$ and ten paths
  $$xy_{i,1}^1y_{i,1}^2...y_{i,1}^{(11t-1)/2}y_i^{(31t-115)/2+i}$$
  $$xy_{i,2}^1y_{i,2}^2...y_{i,2}^{(13t-1)/2}y_i^{(51t-103)/2+2i}$$
  $$xy_{i,3}^1y_{i,3}^2...y_{i,3}^{(13t-1)/2}y_i^{(71t+315)/2+3i}$$
  $$xy_{i,4}^1y_{i,4}^2...y_{i,4}^{(15t-1)/2}y_i^{(91t+313)/2+4i}$$
  $$xy_{i,5}^1y_{i,5}^2...y_{i,5}^{(15t-1)/2}y_i^{(111t+313)/2+5i}$$
  $$xy_{i,6}^1y_{i,6}^2...y_{i,6}^{(17t-1)/2}y_i^{(131t+311)/2+6i}$$
  $$xy_{i,7}^1y_{i,7}^2...y_{i,7}^{(17t-1)/2}y_i^{(151t+309)/2+7i}$$
  $$xy_{i,8}^1y_{i,8}^2...y_{i,8}^{(19t-1)/2}y_i^{(171t+297)/2+8i}$$
  $$xy_{i,9}^1y_{i,9}^2...y_{i,9}^{(19t-1)/2}y_i^{(191t+301)/2+9i}$$
  $$xy_{i,10}^1y_{i,10}^2...y_{i,10}^{(21t-1)/2}y_i^{(211t+305)/2+10i}.$$

\par
Based on the construction,  $B_{21t+i-57}$ contains exactly
sixty-six cycles of lengths:
 \par
 $$\begin{array}{llll}21t+i-57,& 22t+i+7,& 23t+i+210,& 24t+i,\\
  25t+i+1,& 26t+i,& 27t+i,& 28t+i-5,\\
  29t+i+3,& 30t+i+3,& 31t+i+742,& 32t+2i-51,\\
  32t+2i+216,& 34t+2i+209,& 34t+2i,& 36t+2i,\\
  36t+2i-1,& 38t+2i-6,& 38t+2i-3,& 40t+2i+5,\\
  40t+2i+744,& 42t+3i+158,& 43t+3i+215,& 44t+3i+209,\\
  45t+3i-1,& 46t+3i-1,& 47t+3i-7,& 48t+3i-4,\\
  49t+3i-1,& 50t+3i+746,& 53t+4i+157,& 53t+4i+215,\\
  55t+4i+208,& 55t+4i-2,& 57t+4i-7,& 57t+4i-5,\\
  59t+4i-2,& 59t+4i+740,& 63t+5i+157,& 64t+5i+214,\\
  65t+5i+207,& 66t+5i-8,& 67t+5i-5,& 68t+5i-3,\\
  69t+5i+739,& 74t+6i+156,& 74t+6i+213,& 76t+6i+201,\\
  76t+6i-6,& 78t+6i-3,& 78t+6i+738,& 84t+7i+155,\\
  85t+7i+207,& 86t+7i+203,& 87t+7i-4,& 88t+7i+738,\\
  95t+8i+149,& 95t+8i+209,& 97t+8i+205,& 97t+8i+737,\\
  105t+9i+151,& 106t+9i+211,& 107t+9i+946,& 116t+10i+153,\\
  116t+10i+952,& 126t+11i+894.&&\end{array}$$
 \par
 $B_{0}$ is a path with an end vertex $x$ and length $n-n_{t}$. Other $B_i$ is
 simply a cycle of length $i$.
 \par
 It is easy to see that
 $$\begin{array}{lll} v(G)&=& v(B_0)+\sum_{i=1}^{19t+\frac{7t+1}{4}}(v(B_i)-1)+
 \sum_{i=\frac{7t+1}{8}}^{t-742}(v(B_{19t+2i+1})-1)\\

&&+\sum_{i=\frac{7t+1}{8}}^{t-742}(v(B_{19t+2i+2})-1)+\sum_{i=21t-1481}^{21t}(v(B_i)-1)\\
 &&+
 \sum_{i=58}^\frac{7t-7}{8}(v(B_{21t+i-57})-1)
 +\sum_{i=22t-798}^{22t+64}(v(B_i)-1)
 +\sum_{i=23t-734}^{23t+267}(v(B_i)-1)\\
 &&+\sum_{i=24t-531}^{24t+57}(v(B_i)-1)
 +\sum_{i=25t-741}^{25t+58}(v(B_i)-1)
 +\sum_{i=26t-740}^{26t+57}(v(B_i)-1)\\
 &&+\sum_{i=27t-741}^{27t+57}(v(B_i)-1)
 +\sum_{i=28t-741}^{28t+52}(v(B_i)-1)
 +\sum_{i=29t-746}^{29t+60}(v(B_i)-1)\\
 &&+\sum_{i=30t-738}^{30t+60}(v(B_i)-1)
 +\sum_{i=31t-738}^{31t+799}(v(B_i)-1)\end{array}$$
 $$\begin{array}{lll}&=& n-n_t+1+\sum_{i=1}^{19t+\frac{7t+1}{4}}(i-1)+
 \sum_{i=\frac{7t+1}{8}}^{t-742}(144t+13i+1463\\
 &&+\frac{11t-1}{2}+\frac{13t-1}{2}+\frac{13t-1}{2}+
 \frac{15t-1}{2}+\frac{15t-1}{2}+\frac{17t-1}{2}+\frac{17t-1}{2}\\
 &&+\frac{19t-1}{2}+\frac{19t-1}{2}+\frac{21t-1}{2}+\frac{t-571}{2})
+\sum_{i=\frac{7t+1}{8}}^{t-742}(19t+2i+1)\\
&&+\sum_{i=21t-1481}^{21t}(i-1)+\sum_{i=58}^\frac{7t-7}{8}(126t+11i+893\\
&&+\frac{11t-1}{2}+\frac{13t-1}{2}+\frac{13t-1}{2}+
 \frac{15t-1}{2}+\frac{15t-1}{2}+\frac{17t-1}{2}+\frac{17t-1}{2}\\
 &&+\frac{19t-1}{2}+\frac{19t-1}{2}+\frac{21t-1}{2})
 +\sum_{i=22t-798}^{22t+64}(i-1)\\
 &&+\sum_{i=23t-734}^{23t+267}(i-1)+\sum_{i=24t-531}^{24t+57}(i-1)+
 \sum_{i=25t-741}^{25t+58}(i-1)\\
 &&+\sum_{i=26t-740}^{26t+57}(i-1)+\sum_{i=27t-741}^{27t+57}(i-1)
 +\sum_{i=28t-741}^{28t+52}(i-1)\\
 &&+\sum_{i=29t-746}^{29t+60}(i-1)+\sum_{i=30t-738}^{30t+60}(i-1)
 +\sum_{i=31t-738}^{31t+799}(i-1)\\
 &=&n-n_t+\frac{1}{16}(-3309665+1028882t+6911t^2)\\
 &=&n.\end{array}$$
 \par
 Now we compute the number of edges of $G$
  $$\begin{array}{lll} \epsilon(G)&=&
\epsilon(B_0)+\sum_{i=1}^{19t+\frac{7t+1}{4}}\epsilon(B_i)+
 \sum_{i=\frac{7t+1}{8}}^{t-742}\epsilon(B_{19t+2i+1})\\

&&+\sum_{i=\frac{7t+1}{8}}^{t-742}\epsilon(B_{19t+2i+2})+\sum_{i=21t-1481}^{21t}\epsilon(B_i)\\
 &&+
 \sum_{i=58}^\frac{7t-7}{8}\epsilon(B_{21t+i-57})
 +\sum_{i=22t-798}^{22t+64}\epsilon(B_i)
 +\sum_{i=23t-734}^{23t+267}\epsilon(B_i)\\
 &&+\sum_{i=24t-531}^{24t+57}\epsilon(B_i)
 +\sum_{i=25t-741}^{25t+58}\epsilon(B_i)
 +\sum_{i=26t-740}^{26t+57}\epsilon(B_i)\\
 &&+\sum_{i=27t-741}^{27t+57}\epsilon(B_i)
 +\sum_{i=28t-741}^{28t+52}\epsilon(B_i)
 +\sum_{i=29t-746}^{29t+60}\epsilon(B_i)\\
 &&+\sum_{i=30t-738}^{30t+60}\epsilon(B_i)
 +\sum_{i=31t-738}^{31t+799}\epsilon(B_i)\\
 &=& n-n_t+\sum_{i=1}^{19t+\frac{7t+1}{4}}i+
 \sum_{i=\frac{7t+1}{8}}^{t-742}(144t+13i+1464\\
 &&+\frac{11t+1}{2}+\frac{13t+1}{2}+\frac{13t+1}{2}+
 \frac{15t+1}{2}+\frac{15t+1}{2}+\frac{17t+1}{2}+\frac{17t+1}{2}\\
 &&+\frac{19t+1}{2}+\frac{19t+1}{2}+\frac{21t+1}{2}+\frac{t-571+2}{2})
+\sum_{i=\frac{7t+1}{8}}^{t-742}(19t+2i+2)\\
&&+\sum_{i=21t-1481}^{21t}i+\sum_{i=58}^\frac{7t-7}{8}(126t+11i+894\\
&&+\frac{11t+1}{2}+\frac{13t+1}{2}+\frac{13t+1}{2}+
 \frac{15t+1}{2}+\frac{15t+1}{2}+\frac{17t+1}{2}+\frac{17t+1}{2}\\
 &&+\frac{19t+1}{2}+\frac{19t+1}{2}+\frac{21t+1}{2})
 +\sum_{i=22t-798}^{22t+64}i\\
 &&+\sum_{i=23t-734}^{23t+267}i+\sum_{i=24t-531}^{24t+57}i+
 \sum_{i=25t-741}^{25t+58}i\\
 &&+\sum_{i=26t-740}^{26t+57}i+\sum_{i=27t-741}^{27t+57}i
 +\sum_{i=28t-741}^{28t+52}i\\
 &&+\sum_{i=29t-746}^{29t+60}i+\sum_{i=30t-738}^{30t+60}i
 +\sum_{i=31t-738}^{31t+799}i\\
 &=&n-n_t+\frac{1}{16}(-3309681+1029394t+6911t^2)\\
 &=&n+32t-1.\end{array}$$
  \par Then  $f(n)\geq n+32t-1,$ for $n\geq n_{t}.$ This completes
  the proof of the theorem.
 \vskip 0.2in
 \par
  From the above theorem, we have $$\liminf_{n \rightarrow \infty} {f(n)-n\over
\sqrt n}
  \geq \sqrt{2+{2562\over 6911}},$$
  which is better than the previous bounds $\sqrt 2$ (see [2]), $\sqrt
{2+{487\over 1381}}$
  (see [6]).
  \par
  Combining this with Boros, Caro, F\"uredi and Yuster's upper bound, we have
  $$1.98\geq \limsup_{n \rightarrow \infty} {f(n)-n\over \sqrt n} \geq
\liminf_{n \rightarrow \infty} {f(n)-n\over \sqrt n}\geq 1.5397.$$
 \section*{Acknowledgment}
 The author thanks Prof. Yair Caro and Raphael Yuster for sending reference [7].
 The author also thanks Prof. Cheng Zhao for his advice.
 \par

\end{document}